\documentclass[12pt]{amsart}
\usepackage{amssymb}
\usepackage{epsf, epsfig}
\usepackage{multirow}
\usepackage{colordvi}
\numberwithin{equation}{section}

\newtheorem{thm}[equation]{Theorem}
\newtheorem{prop}[equation]{Proposition}

\newtheorem{conj}[equation]{Conjecture}
\newtheorem{example}[equation]{Example}
\newtheorem{remark}[equation]{Remark}
\newtheorem{definition}[equation]{Definition}

\newenvironment{rem}{\begin{remark}\rm}{\end{remark}}
\newcounter{FNC}[page]
\def\newfootnote#1{{\addtocounter{FNC}{2}$^\fnsymbol{FNC}$%
     \let\thefootnote\relax\footnotetext{$^\fnsymbol{FNC}$#1}}}

\newcommand{\Flan}{{{\mathbb F}\ell(a;n)}}
\newcommand{\Fl}{{{\mathbb F}\ell}}
\newcommand{\Gr}{{\rm Gr}}

\newcommand{\C}{{{\mathbb C}}}
\renewcommand{\P}{{{\mathbb P}}}
\newcommand{\R}{{{\mathbb R}}}
\newcommand{\RP}{{\R\P}}

\newcommand{\calM}{{{\mathcal M}}}

\newcommand{\Edot}{{E_\bullet}}
\newcommand{\Fdot}{{F_\bullet}}

\newcommand{\Lamw}[1]{{\mathcal L}_{a_{#1},m}(w)}

\copyrightinfo{}{}
\title[Experimentation in the real Schubert
 calculus]{Experimentation and Conjectures in the real Schubert
 calculus for flag manifolds}

\author{Jim Ruffo}
\address{Department of Mathematics\\
         Texas A\&M University\\
         College Station\\
         TX \ 77843\\
         USA}
\email{ruffo@math.tamu.edu}
\urladdr{http://www.math.tamu.edu/\~{}ruffo}
\author{Yuval Sivan}
\address{Department of Mathematics\\
        University of Massachusetts\\
        Amherst, MA, 01003\\
        USA}
\email{yuval@student.umass.edu}
\author{Evgenia Soprunova}
\address{Department of Mathematics\\
        University of Massachusetts\\
        Amherst, MA, 01003\\
        USA}
\email{soproun@math.umass.edu}
\urladdr{http://www.math.umass.edu/\~{}soproun}
\author{Frank Sottile}
\address{Department of Mathematics\\
         Texas A\&M University\\
         College Station\\
         TX \ 77843\\
         USA}
\email{sottile@math.tamu.edu}
\urladdr{http://www.math.tamu.edu/\~{}sottile}

\thanks{Work done at MSRI and supported by NSF grant DMS-9810361}
\thanks{Computations done on computers purchased with NSF SCREMS grant DMS-0079536}
\thanks{Work of Sottile was supported by the Clay Mathematical Institute}
\thanks{This work was supported in part by NSF CAREER grant DMS-0134860}

\begin{document}
\begin{abstract}
 The Shapiro conjecture in the real Schubert calculus fails to hold for flag
 manifolds, but in a very interesting way.
 We give a refinement of that conjecture for the flag manifold and present
 massive experimentation that supports this conjecture.
 We also establish
 relationships between different cases of the conjecture.
\end{abstract}
\maketitle


\section*{Introduction}
 The Shapiro conjecture for Grassmannians~\cite{So00a,SS02} has driven
 progress in enumerative real algebraic geometry~\cite{So03b}, which is the
 study of real solutions to geometric problems. 
 It conjectures that a (zero-dimensional)
 intersection of Schubert subvarieties of a Grassmannian consists
 entirely of real points---if the Schubert subvarieties are given by
 flags osculating a real rational normal curve.
 This particular geometric problem can be interpreted in terms of real linear series on
 $\P^1$ with  prescribed (real) ramification~\cite{EH83,EH87}, real rational
 curves in $\P^n$ with real flexes~\cite{KhS03}, linear systems
 theory~\cite{RS98}, and the Bethe Ansatz and Fuchsian equations~\cite{MV04}. 
 The Shapiro conjecture has implications for all these areas.
 Massive computational evidence~\cite{So00a,Ver00} as well as its
 proof by Eremenko and Gabrielov for Grassmannians of codimension 2
 subspaces~\cite{EG02a} give compelling evidence
 for its validity. 
 A local version, that it holds when the Schubert varieties are 
 special (a technical term) and when the points of osculation are
 sufficiently clustered~\cite{So99a}, showed that the special
 Schubert calculus is fully real (such geometric problems can have
 all their solutions real).
 That result was generalized by Vakil~\cite{Va04} (using other
 methods) to show that the general Schubert calculus on the
 Grassmannian is fully real.

 The original Shapiro conjecture was that such an intersection of Schubert
 varieties in a \Blue{\emph{flag manifold}} would consist entirely of real
 points.
 Unfortunately, this fails for the first non-trivial
 enumerative problem on a non-Grassmannian flag manifold (see
 Section~\ref{S:F124}), but in a very interesting way. 
 This failure was first noted in~\cite[\S 5]{So00a} and a more symmetric
 counterexample was found in~\cite{So00b}, where computer experimentation 
 suggested that the conjecture would hold if the points where the flags
 osculated the rational normal curve satisfied a certain non-crossing condition.
 Further experimentation led to a precise formulation of this refined 
 non-crossing conjecture in~\cite{So03b}.
 That conjecture was only valid for two- and three- step flag manifolds, and the further
 experimentation reported here leads to a version (Conjecture~\ref{C:Main}) for all flag
 manifolds in which the points of osculation satisfy a monotonicity property.

 We systematically investigate the Shapiro conjecture for flag manifolds to gain
 a deeper understanding both of its failure and of our refinement.
 Our investigation includes 12 gigaHertz-years of computer experimentation and
 theorems relating our monotone conjecture for different enumerative problems.
 Recently, our conjecture was proven by Eremenko, Gabrielov,
 Shapiro, and Vainshtein~\cite{EGSV} for manifolds of flags
 consisting of a codimension 2 plane lying on a hyperplane.

 Our conjecture is concerned with a subclass of Schubert intersection problems.
 Here is one open instance of this conjecture, expressed as a
 system of polynomials in local coordinates for the variety of flags
 $E_2\subset E_3$ in $5$-space, where $\dim E_i=i$.
 Let $t,x_1,\dotsc,x_8$ be indeterminates, and consider the polynomials
\[
   f(t;x)\ :=\   \det\left[
   \begin{array}{ccccc}
    1 & 0 & x_1 & x_2 & x_3\\
    0 & 1 & x_4 & x_5 & x_6\\
    t^4 & t^3 & t^2 & t & 1\rule{0pt}{14pt}\\
    4t^3 & 3t^2 & 2t &1 & 0\\
    6t^2 & 3t &  1  & 0 & 0
   \end{array}\right]\,,
    \qquad
   g(t;x)\ :=\   \det\left[
   \begin{array}{ccccc}
    1 & 0 & x_1 & x_2 & x_3\\
    0 & 1 & x_4 & x_5 & x_6\\
    0 & 0 &  1  & x_7 & x_8\\
    t^4 & t^3 & t^2 & t & 1\rule{0pt}{14pt}\\
    4t^3 & 3t^2 & 2t &1 & 0\\
   \end{array}\right]\, .
\]

\noindent{\bf Conjecture A.}
{\it
  Let $t_1<t_2<\dotsb<t_8$ be real numbers.
  Then the polynomial system
\begin{eqnarray*}
  \Blue{f(t_1;x)=f(t_2;x)=f(t_3;x)=f(t_4;x)}& =& 0,\quad\mbox{and}\\ 
  \Red{g(t_5;x)=g(t_6;x)=g(t_7;x)=g(t_8;x)}&=& 0
\end{eqnarray*}
 has $12$ solutions, and all of them are real.
}\medskip

In Conjecture A, monotonicity is that the polynomials $f$ are evaluated at parameter
values that are less than the parameter values at which the polynomials
$g$ are evaluated.
If the order of $t_4$ and $t_5$ were switched, then the evaluation
would not be monotone.
We computed 100,000 instances of this polynomial system at different
ordered parameter values, and each had 12 real solutions.
In contrast, we found non-monotone evaluations for which not all solutions were 
real, and the minimum number of real solutions that we observed
depended on the combinatorics of the evaluation. 
This is summarized in Table~\ref{table:12-flag}.\smallskip

Section 1 contains background on flag manifolds, states the Shapiro Conjecture, and
gives a geometrically vivid example of its failure.
In Section 2, we state our conjectures and discuss relations among them.
The discussion in Section 3 contains theorems about the conjectures.
Finally, in Section 4 we describe our methods, explaining our experimentation and giving
a brief guide to our data, all of which and much more is tabulated and available online
at {\tt www.math.tamu.edu/\~{}sottile/pages/Flags/}.
We also describe some interesting phenomena we observed in our data.

\section{Background}

\subsection{Basics on flag manifolds}
Given positive integers $a:=(a_1<\dotsb<a_k)$ with $a_k<n$, let $\Flan$ 
be the manifold of flags in $\C^n$ of {\it type} $a$,
\[
   \Flan\ :=\ \{\Edot=E_{a_1}\subset E_{a_2}\subset\dotsb\subset
                 E_{a_k}\subset\C^n\mid \dim E_{a_i}=a_i\}\,.
\]
If we set $a_0:=0$, then this has dimension 
$ \dim(a) :=  \sum_{i=1}^k (n-a_i)(a_i-a_{i-1})$.
Flags $\Fdot$ of type $1<2<\dotsb<n{-}1$ in $\C^n$ are called {\it complete}.

The positions of flags $\Edot$ of type $a$ relative to a fixed
complete flag $\Fdot$ stratify $\Flan$ into {\em Schubert cells} whose closures are 
{\em Schubert varieties}.
These have a precise description in terms of linear
algebra and combinatorics.
Define $W^a\subset S_n$ to be the set of permutations of $\{1,2,\dotsc,n\}$ with descents
in $a$,
\[ 
   W^a\ :=\ \{w\in S_n\mid i\not\in\{a_1,\dotsc,a_k\}\Rightarrow  
      w(i)<w(i+1)\}\,.
\]
Permutations $w\in W^a$ index Schubert cells $X^\circ_w\Fdot$ and 
Schubert varieties $X_w\Fdot$.
Precisely, if we set $r_w(i,j)=|\{l\leq i\mid j+w(l)>n\}|$, then 
 \begin{eqnarray}
   X^{\Red{\circ}}_w\Fdot&=&\{\Edot\mid \dim E_{a_i}\cap F_j\Red{=}\nonumber 
       r_w(a_i,j),\ i=1,\dotsc,k,\ j=1,\dotsc,n\},\quad\mbox{and}\\ 
  X_w\Fdot&=&\{\Edot\mid \dim E_{a_i}\cap F_j\Red{\geq} r_w(a_i,j),\
       i=1,\dotsc,k,\ j=1,\dotsc,n\}\,.
 \end{eqnarray}
Flags $\Edot$ in $X^\circ_w\Fdot$ have {\em position $w$ relative to $\Fdot$}.
We call a permutation $w\in W^a$ a {\em Schubert condition} on flags of type
$a$. 
These irreducible subvarieties have codimension in $\Flan$ equal to the
length, $\ell(w)$, of the index $w$.

We have $X^\circ_w\Fdot\simeq \C^{\dim(a)-\ell(w)}$.
We use a convenient set of coordinates for the Schubert cells.
Let $\calM_w$ be the set of  $a_k\times n$-matrices $M$ whose
entries $x_{i,j}$ satisfy
 \begin{eqnarray*}
   x_{i,w(i)}& =& 1\\
   x_{i,j}& =& \makebox[30pt][l]{0}\mbox{if}\ j<w(i)\ \mbox{or}\ w^{-1}(j)<i\\
       &&\makebox[30pt][r]{ or\ } \mbox{if } a_l<i<w^{-1}(j)\leq a_{l+1},\ 
            \mbox{for some }l\,,
 \end{eqnarray*}
and whose remaining $\dim(a)-\ell(w)$ entries 
give coordinates for $\calM_w$.
For example, if $n=8$, $a=(2,3,6)$, and $w=25\,3\,167\,48$, then 
$\calM_w$ consists of matrices of the form
\[
  \left(\begin{matrix}
      0   &1&x_{13}&x_{14}&0&x_{16}&x_{17}&x_{18}\\
      0   &0&  0   &   0  &1&x_{26}&x_{27}&x_{28}\\
      0   &0&  1   &x_{34}&0&x_{36}&x_{37}&x_{38}\rule{0pt}{13pt}\\
      1   &0&  0   &x_{44}&0&  0   &  0   &x_{48}\rule{0pt}{13pt}\\
      0   &0&  0   &   0  &0&  1   &  0   &x_{58}\\
      0   &0&  0   &   0  &0&  0   &  1   &x_{68}
   \end{matrix}\right)\ .
\]

The relation of $\calM_w$ to the Schubert cell is as follows.
Given a complete flag $\Fdot$, choose an ordered basis 
$e_1,\dotsc,e_n$ for $\C^n$ corresponding to
the columns of our matrices, such that 
$F_i$ is the linear span of the last $i$ basis vectors,
$e_{n+1-i},\dotsc,e_{n-1},e_n$. 
Set $E_{a_i}$ to be the row space of the
first $a_i$ rows of a matrix $M\in\calM_w$.
Then the flag $\Edot$ has type $a$ and lies in the 
Schubert cell $X^\circ_w\Fdot$,
every flag $\Edot\in X^\circ_w\Fdot$ arises in this way, and the 
association $M\mapsto \Edot$ is an algebraic bijection 
between $\calM_w$ and $X^\circ_w\Fdot$.

Let $\iota$ be the identity permutation.
Then $\calM_{\iota}$ provides local coordinates for $\Flan$ in which the equations for a
Schubert variety $X_w\Fdot$ are easy to describe.
Note that  
\[
   \dim (E_{a_i}\cap F_j)\ \geq\ r\quad 
    \Longleftrightarrow\quad \mbox{rank}(A)\ \leq\ a_i+j-r\,,
\]
where $A$ is the matrix formed by stacking the first $a_i$ rows of
$\calM_{\iota}$ on top of a $j\times n$ 
matrix for $F_j$.
This rank condition is the vanishing of all minors of $A$ of size 
$1{+}a_i{+}j{-}r$.

When $b=(b)$ is a singleton,  $\mathbb{F}\ell(b;n)$ is the Grassmannian
of $b$-planes in $\C^n$, written $\Gr(b,n)$.
Non-identity permutations in $W^b$ have a unique descent at $b$ and are called
\Blue{\emph{Grassmannian}}. 
Our conjecture concerns \Blue{\emph{Grassmannian Schubert varieties}}, which are Schubert
varieties of any flag manifold indexed by Grassmannian permutations.

\subsection{The Shapiro Conjecture}
A list $(w_1,\dotsc,w_m)$ of permutations in $W^a$ is called 
{\em Schubert data} if $\ell(w_{1})+\cdots+\ell(w_{m})=\dim(a)$.
Schubert data index enumerative geometric problems involving
Schubert varieties. 
Specifically, given such Schubert data and complete flags
$\Fdot^1,\dotsc,\Fdot^m$, consider the Schubert intersection
 \begin{equation} \label{SchInt}
   X_{w_{1}}\Fdot^{1} \cap \dotsb  \cap X_{w_{m}}\Fdot^{m}\,.
 \end{equation}
When the flags are in general position, this intersection is zero-dimensional, 
and it equals the intersection of corresponding Schubert cells.
In that case, the intersection ~\eqref{SchInt} consists of the
flags $\Edot$ of type $a$ which have position $w_i$ relative to $\Fdot^i$, for
each $i=1,\dotsc,m$. 
Counting the solutions is a problem in enumerative geometry.

The degree of a zero-dimensional intersection~\eqref{SchInt} does
not depend on the choice of flags and we call this number $d(w_1,\dotsc,w_m)$ the 
{\em degree} of the Schubert data.
When the intersection is transverse, this degree counts the
solutions to the problem of enumerating flags of type $a$ having
positions $w_1,\dotsc,w_m$ relative to the fixed chosen flags.
Transversality is guaranteed if the chosen flags are in general position~\cite{Kl74}.

The Shapiro conjecture concerns the following variant to this
classical enumerative geometric problem:
What {\em real} flags $\Edot$ have given position $w_i$ relative to
{\em real} flags $\Fdot^i$, for each $i=1,\dotsc,m$?
In the Shapiro conjecture, the flags $\Fdot^i$ are not taken to be general real flags,
but rather flags osculating a rational normal curve.
Let $\gamma\colon \C \rightarrow \C^n$ be the rational normal curve, 
$\gamma\colon t \mapsto (t^{n-1},t^{n-2},\ldots,t,1)$.
The {\em osculating flag} $\Fdot(t)$ of subspaces to the rational
normal curve at the point $\gamma(t)$ is
\[
   \Fdot(t)\ :=\ \mbox{span}\{\gamma(t),
    \frac{d\gamma}{ds}(t),\ldots,\frac{d^{i-1}\gamma}{ds^{i-1}}(t)\}\,.
\]
When $t=0$, the flag $\Fdot(0)$ is the standard
flag we used in the description of $\calM_w$.

\begin{conj}[B.~Shapiro and M.~Shapiro]\label{C:SC}
  Suppose that $w_1,\dotsc,w_m$ is Schubert data for flags of type
  $a$.
  If the flags $\Fdot^1,\dotsc,\Fdot^m$ osculate the
  rational normal curve at distinct real points, then the 
  intersection~$\ref{SchInt}$ is transverse and consists only of real points. 
\end{conj}

Write $X_w(t)$ for $X_w\Fdot(t)$ and $X_w$ for $X_w\Fdot(0)$.
Then the Shapiro conjecture is concerned with intersections of the form
 \begin{equation}\label{E:shapint}
     X_{w_1}(t_1)\cap X_{w_2}(t_2)\cap 
    \dotsb\cap X_{w_m}(t_m)\,.
 \end{equation}

Previously, Conjecture~\ref{C:SC} has mostly been studied for Grassmannians.
Experimental evidence for its validity was first found by one of us
(Sottile)~\cite{RS98,So97c}.
This led to a more systematic investigation, both experimentally and
theoretically~\cite{So00a}.
There, relationships between the conjecture for different collections of
Schubert data on different Grassmannians were studied.
For example, if the Shapiro conjecture holds for a Grassmannian for
the Schubert data consisting only of codimension 1 conditions, then it holds
for all Schubert data on that Grassmannian and on all smaller
Grassmannians.
(Perhaps dropping the claim of transversality.)
More recently, Eremenko and Gabrielov proved the conjecture for
any  Schubert data on a Grassmannian of codimension
2-planes~\cite{EG02a}. 
This result is appealingly interpreted as a rational function with
real critical points must be real.

The original conjecture was for flag manifolds.
Unfortunately, a counterexample was found in~\cite{So00a}.
Subsequent experimentation refined this counterexample, and
suggested a reformulation of the original conjecture that might be
true.
We study this refined conjecture, and  report on
massive computer experimentation (over 12 GHz-years) studying this conjecture.
A by-product of this experimentation was the discovery of several new and
unusual phenomena, which we describe through examples.
The first is a relatively simple counterexample to the Shapiro conjecture.

\subsection{The Shapiro conjecture is false for flags in space}\label{S:F124} 
Use $a^b$ to indicate that condition $a$ is repeated $b$ times.
Then $\bigl((1324)^3,\,(1243)^2\bigr)$ is Schubert data for flags of
type $(2,3)$ in $\C^4$.
For distinct points $s,t,u,v,w\in\RP^1$, consider the Schubert intersection
 \begin{equation}\label{F234Int}
    X_{1324}(s)\cap X_{1324}(t)\cap X_{1324}(u)\  \cap\ 
   X_{1243}(v)\cap X_{1243}(w)\,,
 \end{equation}
in projective 3-space.
There, a partial flag of type $(2,3)$ consists of a line $\ell$ lying on a plane $H$.
Then $(\ell\subset H)\in X_{1243}(v)$ if the plane $H$ contains the point
$\gamma(v)$, and 
$(\ell\subset H)\in X_{1324}(s)$ if $\ell$ meets the tangent line $\ell(s)$ to $\gamma$
at $\gamma(s)$. 

Suppose that the flag $\ell\subset H$ lies in the
intersection~\eqref{F234Int}.
Then $H$ contains the two points $\gamma(v)$ and $\gamma(w)$, and hence the
secant line $\lambda(v,w)$ that they span.
Since $\ell$ is another line in $H$, $\ell$ meets this secant line
$\lambda(v,w)$.
As long as $\ell\neq\lambda(v,w)$, then $\ell$ determines $H$ uniquely as the
span of $\ell$ and $\lambda(v,w)$.
In this way, we are reduced to determining the lines $\ell$ which meet the
three tangent lines $\ell(s)$, $\ell(t)$, $\ell(u)$, and the secant line
$\lambda(v,w)$.

\begin{figure}[htb]
 \[
  \begin{picture}(215,99)(0,5)
   \put(0,0){\includegraphics[height=110pt]{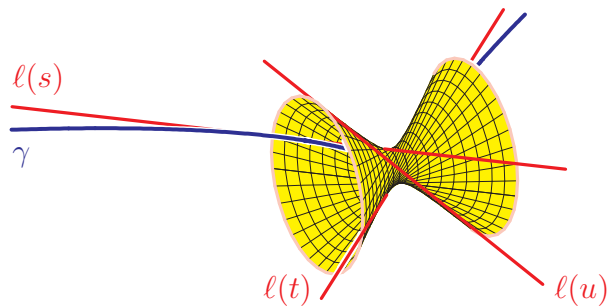}}
   \put(0,80){\Red{$\ell(s)$}} \put(95,0){\Red{$\ell(t)$}} \put(205,0){\Red{$\ell(u)$}}
   \put(0,52){\Blue{$\gamma$}}
  \end{picture}
 \]
 \caption{Quadric containing three lines tangent to the rational normal
 curve\label{F:3TanQuad}} 
 \end{figure}

\begin{figure}[htb]
 \[
  \begin{picture}(200,85)(0,8)
   \put(0,0){\includegraphics[height=100pt]{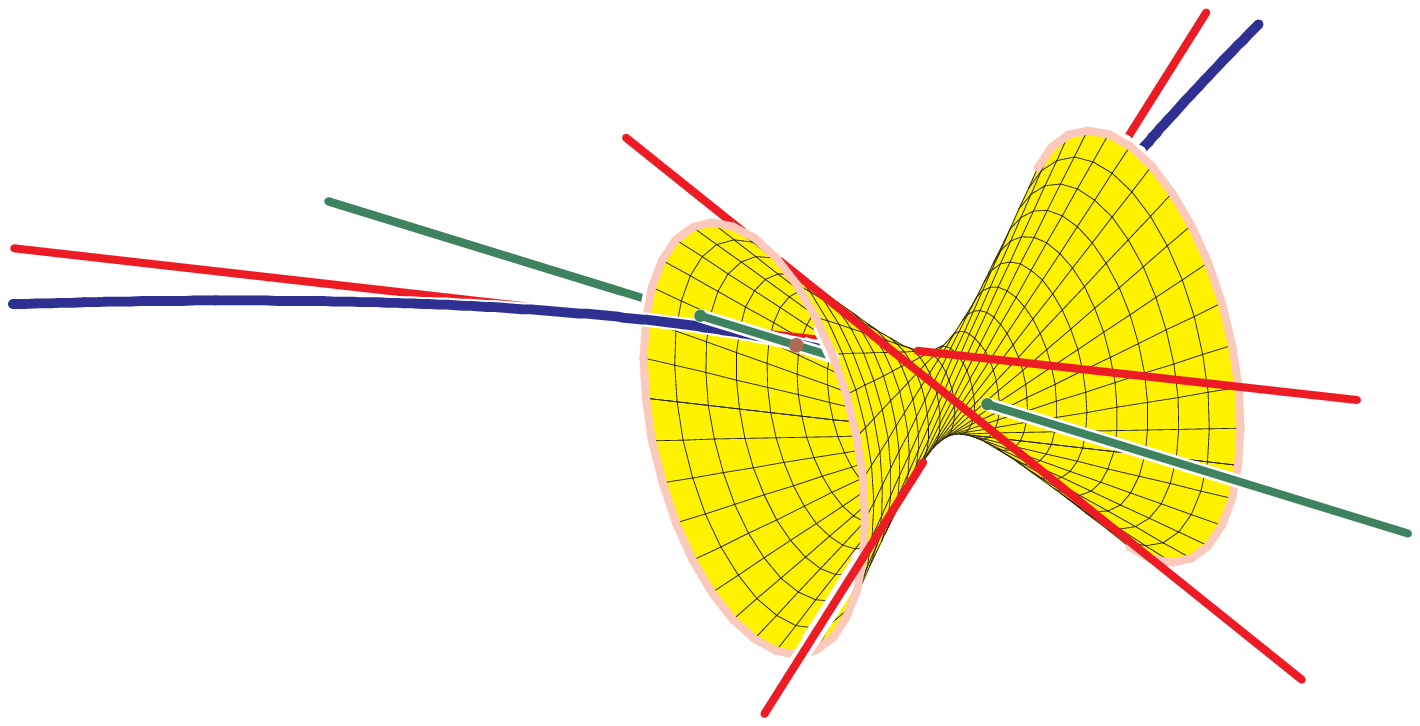}}
   \put(0,72){\Red{$\ell(s)$}} \put(85,0){\Red{$\ell(t)$}} \put(187,0){\Red{$\ell(u)$}}
   \put(0,45){\Blue{$\gamma$}} \put(30,78){\ForestGreen{$\lambda(v,w)$}}
  \end{picture}\qquad
  \begin{picture}(165,85)(-20,5)
   \put(0,0){\includegraphics[height=75pt]{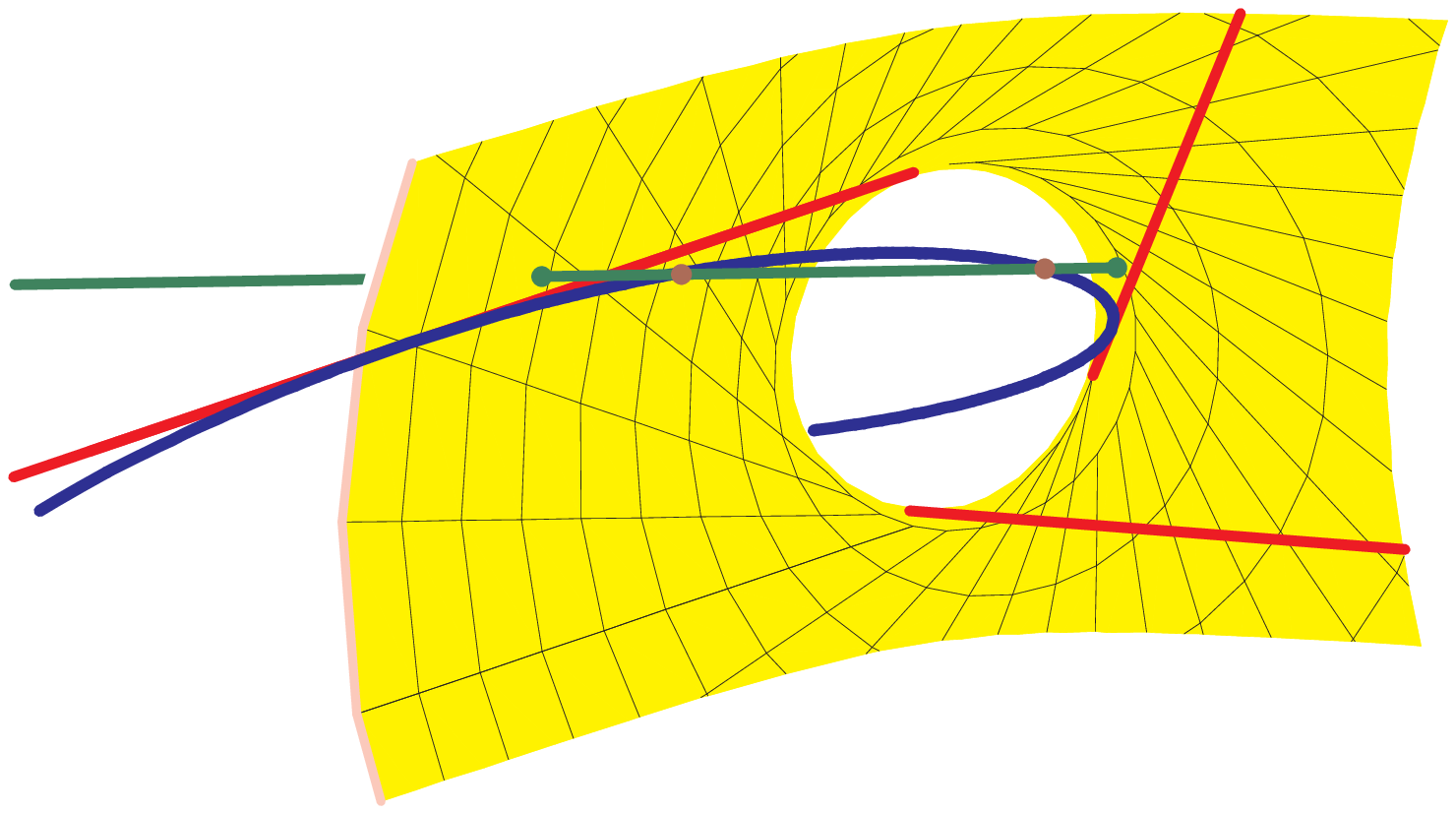} }
   \put(-22,32){\Red{$\ell(s)$}} \put(137,20){\Red{$\ell(t)$}} \put(102,82){\Red{$\ell(u)$}}
   \put(0,20){\Blue{$\gamma$}} \put(-20,57){\ForestGreen{$\lambda(v,w)$}}
  \end{picture}
\] \caption{Two views of secant line meeting $Q$\label{F:NC}}
 \end{figure}

\begin{figure}[htb]
 \[
  \begin{picture}(240,110)(0,10)
   \put(0,0){\includegraphics[height=120pt]{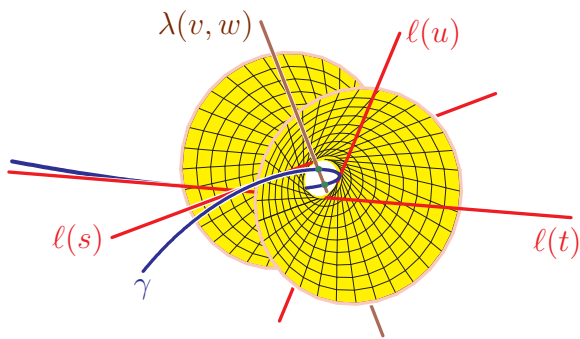} }
   \put(18,34){\Red{$\ell(s)$}} \put(200,33){\Red{$\ell(t)$}} \put(152,112){\Red{$\ell(u)$}}
   \put(49,17){\Blue{$\gamma$}} \put(58,115){\Brown{$\lambda(v,w)$}}
  \end{picture}
 \]
 \caption{Secant line not meeting $Q$\label{F:CR}}
 \end{figure}

The set of lines which meet the three mutually skew lines 
$\ell(s)$, $\ell(t)$, and $\ell(u)$ forms one ruling of a quadric surface $Q$
in $\P^3$. 
We display this quadric $Q$ and the ruling in Figure~\ref{F:3TanQuad}, as well
as the rational normal curve with its three tangent lines. 
The lines meeting these three and the secant line $\lambda(v,w)$ correspond to
the points where $\lambda(v,w)$ meets the quadric $Q$.
In Figure~\ref{F:NC}, we display a secant line which meets the hyperboloid in two
points, and therefore these choices give two real flags in the
intersection~\eqref{F234Int}.
There is also a secant line which meets the hyperboloid in two complex conjugate points.
For this secant line, the two flags in the intersection~\eqref{F234Int} are both
complex. 
We show this in Figure~\ref{F:CR}.

Since any three points of $\RP^1$ may be carried to any other three points via a
real M\"obius transformations, \emph{and} these extend to projective automorphisms of
$\P^3$, we have

\begin{prop}\label{P:counter}
  The intersection~\eqref{F234Int} is transversal and consists only of real
  points if there are disjoint intervals $I_2$ and $I_3$ of\/ $\RP^1$ so that
  $s,t,u\in I_2$ and $v,w\in I_3$.
\end{prop}

Proposition~\ref{P:counter} suggests a refinement to the Shapiro conjecture
which may hold. 
We describe that refinement and experimental evidence in its favor.

\section{Results}

Experimentation designed to test hypotheses is a
primary means of inquiry in the natural sciences.
In mathematics we typically use proof and example as our primary
means of inquiry.
Many mathematicians (including the authors) feel that they are 
striving to understand the nature of objects that inhabit a very
real mathematical reality.
For us, experimentation plays an important role in helping to  
formulate conjectures, which are then
studied and perhaps eventually decided.
We discuss the conjectures which were informed by our massive
experimentation that is described in Section~\ref{S:Methods}.

\subsection{Conjectures}

Let $a=(a_1<\dotsb<a_k)$ and $n$ be positive integers with $a_k<n$.
A permutation $w\in W^a$ is \Blue{\emph{Grassmannian}} if it has a single
descent, say at position $a_l$.
Then the Schubert variety $X_w\Fdot$ of $\Flan$ is the pull-back of the
Schubert variety $X_w\Fdot$ of the Grassmannian $\Gr(a_l,n)$.
Write $\delta(w)$ for the unique descent of  $w$.

Suppose that $w_1,\dotsc,w_m\in W^a$ are data for a Schubert problem on $\Flan$ with  
each permutation $w_i$ Grassmannian.
A collection of points $t_1,\dotsc,t_m\in\R\P^1$ is \Blue{\emph{monotone}}
(with respect to $w_1,\dotsc,w_m$) if the function
\[
   t_i\ \longmapsto\ \delta(w_i)\ \in\{a_1,a_2,\dotsc,a_k\}
\]
is monotone, when the ordering of the
$t_i$ is consistent with some orientation of $\R\P^1$.

\begin{conj}\label{C:Main}
 Suppose that $w_1,\dotsc,w_m\in W^a$ are data for a Schubert problem where 
 each permutation $w_i$ is Grassmannian.
 Then the intersection 
\begin{equation}\label{AA}
     X_{w_1}(t_1)\cap X_{w_2}(t_2)\cap 
    \dotsb\cap X_{w_m}(t_m)\,,
\end{equation}
 is transverse with all points of intersection real, if the points
 $t_1,\dotsc,t_m\in\R\P^1$ are monotone with respect to $w_1,\dotsc,w_m$.
\end{conj}

Note that if $\Flan$ is a Grassmannian, then every choice of points is monotone, so
the Shapiro conjecture for Grassmannians is a special case of Conjecture~\ref{C:Main}.

Our experimentation systematically investigated the original Shapiro conjecture for flag
manifolds, with a focus on this monotone conjecture.
In all, we examined 591 such Grassmannian Schubert problems on 29 different flag
manifolds.
On these, we verified that on each of 140 million specific monotone intersections of the
form~\eqref{AA} had all solutions real.
We find this to be overwhelming evidence in support of our monotone conjecture.

\begin{rem}
 Conjecture~A of the Introduction is the instance of Conjecture~\ref{C:Main} for the
 flag manifold $\mathbb{F}\ell(2<3;5)$ for the Schubert data
 $\bigl( (13245)^4,\, (12435)^4\bigr)$, expressed in terms of polynomial systems in
 local coordinates.
\end{rem}

\begin{rem}
 The example of Section~\ref{S:F124} illustrates both Conjecture~\ref{C:Main} and its
 limitations.
 The condition on disjoint intervals $I_2$ and $I_3$ of Proposition~\ref{P:counter} is
 equivalent to the monotone choice of points in Conjecture~\ref{C:Main}.
 Also note that the choices which give no real solutions are not monotone.
\end{rem}

We give a stronger conjecture which is supported by our
experimental investigation.
It ignores the issue of reality and concentrates only on the 
transversality of an intersection.

\begin{conj}\label{C:transverse}
 If $w_1,\dotsc,w_m\in W^a$ are data for a Schubert problem where 
 each permutation $w_i$ is Grassmannian
 then the intersection 
\[
     X_{w_1}(t_1)\cap X_{w_2}(t_2)\cap 
    \dotsb\cap X_{w_m}(t_m)\,,
\]
 is transverse, if the points
 $t_1,\dotsc,t_m\in\R\P^1$ are monotone with respect to $w_1,\dotsc,w_m$.
\end{conj}

\begin{thm}\label{T:Implies}
 Conjecture~$\ref{C:transverse}$ implies Conjecture~$\ref{C:Main}$.
\end{thm}

We prove this in Section~\ref{S:Discussion}.
In our experimentation, we kept track of which polynomial systems corresponded to
non-transverse intersections.
No such polynomial system came from a monotone choice.

Eremenko,  Gabrielov, Shapiro, and Vainshtein recently established
Conjecture~\ref{C:Main} for manifolds of flags of type $(n{-}2,n{-1})$ in
$\C^n$~\cite{EGSV}.

\section{Discussion}\label{S:Discussion}

For $m\leq n-a_i$, the \Blue{\emph{special}} Schubert 
variety in $\Gr(a_i,n)$ is defined by
\[
 \Omega_{(m)}\Fdot\ :=\ \{E_{a_i}|{\rm dim}~E_{a_i}\cap F_{m}\geq 1\}\,.
\]
The inverse image under the projection map 
$\pi_{a_i}:\Flan\rightarrow \Gr(a_i,n)$
sending $\Edot\mapsto E_{a_i}$ is a Schubert variety
of $\Flan$ defined by the same formula and we use $\Omega_{(m)}\Fdot$ for this
subvariety  as well. 
When $\Fdot=\Fdot(t)$ osculates the rational normal curve at a point
$t\in\P^1$, we denote the corresponding Schubert variety by
$\Omega_{(m)}(t):=\Omega_{(m)}\Fdot(t)$.

$W^{a_i}$ coincides with the set $\binom{[n]}{a_i}$ of subsets of $\{1,\dotsc,n\}$ of
size $a_i$.
For each $i=1,\dotsc,k$, define a map
$\pi_{a_i}:W^a\rightarrow\binom{[n]}{a_i}$ by
$\pi_{a_i}(v):=\{v_1,\dotsc,v_{a_i}\}$.

Define the \Blue{\emph{$a_i$-Bruhat order}} on $W^a$ by its covers:
if $u\lessdot v$ is a cover in the Bruhat order, then 
$u\lessdot_{a_i}v$ if and only if $\pi_{a_i}(u)<\pi_{a_i}(v)$.
Equivalently,  $u\lessdot_{a_i}v$ if and only if  
$v=u \sigma_{bc}$, where $\sigma_{bc}$ is the transposition
interchanging $b\leq a_i<c$.
Note in particular that $u<_{a_i}v$ implies $\pi_{a_i}(u)<\pi_{a_i}(v)$.
Define a set $\Lamw{i}\subset W^a$ by 
 \begin{eqnarray*}
 \Lamw{i}&:=&\{v\in W^a\mid
 w<_{a_i}v,~\ell(v)=\ell(w)+m,~\#\{l>a_i|v(l)\neq w(l)\}=m\}.
\end{eqnarray*}
The following result is established in~\cite{So96}.

\begin{thm}\label{flPieri} 
 In the Chow ring, $A^*(\Flan)$, the cycle-theoretic equality holds:
 \[
 [X_w\Fdot]\cdot[\Omega_{(m)}(t)]=\sum_{v\in\Lamw{i}}[X_v\Fdot]
 \]
 \end{thm}
 
Eisenbud and Harris~\cite[Theorem 8.1]{EH83}
proved that the scheme-theoretic limit
$\lim_{s\to 0} X_w\Fdot\cap \Omega_{(m)}(t)$ 
is supported on a union of Schubert
varieties, when the \emph{flag manifold is a Grassmannian}.
We remove this restriction to Grassmannians.

\begin{thm}\label{limthm}
 The equality
 \begin{equation}\label{lim}
  \lim_{t\rightarrow 0}~X_w\Fdot\cap\Omega_{(m)}(t)=\sum_{v\in\Lamw{i}}X_v\Fdot
 \end{equation}
 holds scheme-theoretically.
\end{thm}

Our proof follows that of~\cite{EH83}. 
This was proven when $m=1$ in~\cite{So00b}.
It has a number of consequences for the monotone conjecture.

\begin{thm}\label{T:Limits} 
 If Conjecture~$\ref{C:Main}$ holds for all Schubert data on a
 flag manifold $\Flan$ involving codimension $1$ conditions
 (each permutation $w_i$ has length 1), then it holds for all Schubert 
 data on $\Flan$.
\end{thm}

\begin{proof}[Proof of Theorem~$\ref{T:Implies}$]
If Conjecture~\ref{C:transverse} holds, then by~\cite[Corollary 2.2]{So00b},
Conjecture~\ref{C:Main} holds for all Schubert intersection 
problems~\eqref{AA} involving only codimension 1 Schubert conditions.
By Theorem~\ref{T:Limits}, this implies Conjecture~\ref{C:Main} for all Schubert
intersection problems.
\end{proof}

\section{Methods}\label{S:Methods}

At the heart of our experiments is the computation the number of 
real roots of a real polynomial system.  This system generates the
ideal of the Schubert intersection
\[
 X_{w_1}(t_1)\cap X_{w_2}(t_2)\cap\dotsb\cap X_{w_m}(t_m)\,,
\]
which depends upon the points
$t_1,\dotsc,t_m\in\R\P^1$ at which the fixed flags osculate the 
curve $\gamma$.
Our computational procedure takes the following data as input:
\begin{enumerate}
    \item A compact description of the Schubert problem.
    \item The number of complex solutions to the Schubert problem.
    \item A list of \Blue{\emph{necklaces}}, that is, the combinatorial types of the
          order along the rational normal curve at which the conditions will be evaluated.
    \item A list, $L\subset\R\P^1$, of points 
          where the conditions are to be evaluated.
\end{enumerate}
These data are used to create a 
Singular~\cite{SINGULAR} input file that controls 
the first stage of the computation. In this stage, a subset of the 
points $L$ is selected randomly, then each necklace determines which 
Schubert conditions are evaluated at which osculating flags. Each such 
choice is one instance of the Schubert problem, whose ideal is written in the 
Singular input file. 
Singular is called, and it uses Gr\"{o}bner bases to compute 
an eliminant (a univariate polynomial in the ideal of the Schubert
intersection) for each instance. 
Finally, Maple is called to compute 
the number of real roots of each eliminant using the package 
{\tt realroot}, and a table (described below) is updated accordingly.
Instances for which an eliminant could not be computed were
set aside and studied by hand. 

One iteration of the procedure just described yields one instance 
of the Schubert problem for each necklace. The complete computation
is organized by a shell script which iterates this procedure a fixed 
number of times (typically several hundred to several tens of thousands), 
randomly selecting a new set of points $L$ each time.  

Table~\ref{table:12-flag} shows the results of computing 800,000 
 \begin{table}[htb]
  \begin{tabular}
   {|c||c|c|c|c|c|c|c|}\hline
   {Necklace} & \multicolumn{7}{c|}{Number of Real Solutions\rule{0pt}{13pt}}\\
   \cline{2-8}
        &0&2&4&6&8&10&12\rule{0pt}{13pt}\\\hline\hline
   22223333
     &0 & 0 & 0 & 0 & 0 & 0 & 100000\\\hline
   22322333   
     &0 & 0 & 21 & 16129 & 33686 & 29350 & 20814\\ \hline
   22233233
     &0 & 0 & 31 & 16276 & 33430 & 29194 & 21069 \\ \hline
   22332233
     &0 & 0 & 421 & 13742 & 46961 & 23561 & 15315  \\ \hline
   22323323
     &0 & 0 & 6242 & 22480 & 36329 & 26522 & 8427 \\ \hline
   22332323
     &0 & 504 & 18532 & 27844 & 30962 & 15546 & 6612\\ \hline
   22232333
     &0 & 1846 & 8414 & 13887 & 25079 & 20784 & 29990 \\ \hline
   23232323
     & 3830 & 10131 & 32326 & 21679 & 20790 & 8066 & 3178 \\ \hline
  \end{tabular}\vspace{5pt}

  \caption{Table for the problem $13245^4\cdot 12435^4=12$ on $\Fl(2,3;5)$} 
  \label{table:12-flag}
 \end{table}
instances of the Schubert problem $(13245^4,12435^4)$ on 
$\Fl(2 < 3;5)$.  Each row records the number of times a
given number of real solutions was observed for a given necklace.  The
entries in the first column represent the necklaces as sequences 
$\{\delta(w_1),\dotsc,\delta(w_m)\}$, where $\delta(w)$ denotes the unique
descent of the Grassmannian permutation $w$, as described in Section $2.1$.
Thus, in this case, a $2$
represents the condition on the $2$ plane $E_2$ given by the permutation
$1324$, and a $3$ represents the condition on $E_3$ given by the permutation
$1243$.  The necklace corresponding to the first row is monotone.

\subsection{Observations}

This experimentation not only studied Conjecture~\ref{C:Main}, but it systematically
studied the original Shapiro conjecture.
Many interesting phenomena were observed.
For example, there is an extension of Conjecture~\ref{C:Main} in which some conditions
are not Grassmannian (and a similar extension of Theorem~\ref{T:Limits}).
We do not state it here in the extended abstract.

We found many Schubert problems with an apparent \emph{lower bound} on their number of
real solutions.  
For example, Table~\ref{table:lower} is from the problem on $\Fl(2,3;6)$
with Schubert data $(132456^5,125346^3)$, which has degree 14.
 \begin{table}[htb]
  \begin{tabular}
  {|c||c|c|c|c|c|c|c|c|}\hline
   {Necklace} & \multicolumn{8}{c|}{Number of Real Solutions\rule{0pt}{13pt}}\\
   \cline{2-9}
        &0&2&4&6&8&10&12&14\rule{0pt}{13pt}\\\hline\hline
  22222333 &  0  &  0  &   0  &   0  &   0 &   0 &    0 & 3000 \\\hline
  22223233 &  0  &  0  &   0  &   0  &  27 & 562 &  887 & 1524 \\\hline
  22232323 &  0  &  0  &   4  &  77  & 474 & 776 &  810 &  859 \\\hline
  22232233 &  0  &  0  &  34  &  69  & 331 & 915 & 1063 &  588 \\\hline
  22322323 &  0  &  0  & 152  & 356  & 634 & 839 &  726 &  293 \\ \hline
\end{tabular} \vspace{5pt}
  \caption{Table for the problem $132456^5\cdot 125346^3=14$ on $\Fl(2,3;6)$} 
  \label{table:lower}
 \end{table}
Such lower bounds on the number of real solutions to enumerative geometric problems were
first found by Eremenko and Gabrielov~\cite{EG01b} in the context of the Shapiro
conjecture for Grassmannians.
Lower bounds have also been observed for rational curves on surfaces~\cite{IKS,Mi,W}
and for sparse polynomial systems~\cite{SS}.

The problem $(312564^2, 124356^5)$ on $\Fl(1,3,5;6)$ has degree 10, and the condition
$ 312564$ is \emph{not} Grassmannian and $124356$ is Grassmannian with descent at 3.
Not only does this problem exhibit a lower bound, but it also has apparent `gaps' in the
possible numbers of real solutions.
Table~\ref{T:gaps} gives the data from this computation.
In each necklace, $A$ represents the condition $312564$, while 3 represents the
Grassmannian condition $124356$.
 \begin{table}[htb]
  \begin{tabular}
    {|c||c|c|c|c|c|c|}\hline
     {Necklace} & \multicolumn{6}{c|}{Number of Real Solutions\rule{0pt}{13pt}}\\
     \cline{2-7}
        &0&2&4&6&8&10\rule{0pt}{13pt}\\\hline\hline
     $A$33$A$333 & 0 &  1850 & 0 & 10381 & 0 & 7769  \\ \hline
     $A$3$A$3333 & 0 &  3177 & 0 & 13729 & 0 & 3094  \\ \hline 
     $A$$A$33333 & 0 & 11222 & 0 &  8397 & 0 &  381  \\ \hline 
  \end{tabular} \vspace{5pt}
  \caption{Table for the problem $312564^2\cdot 124356^5=10$ on $\Fl(1,3,5;6)$}
  \label{T:gaps}
 \end{table}
This is a new phenomena first observed in some sparse polynomial systems~\cite{SS}.

One unusual problem we looked at was on the flag manifold $\Fl(2,4;6)$
and it involved four identical, but non-Grassmannian conditions,
$142536$.
It always had real solutions, and by using some geometric reasoning we can prove
that is always the case.

\begin{thm}
  For any distinct $s,t,u,v\in\R\P^1$, then intersection
\[ 
   X_{142536}(s) \cap X_{142536}(t) \cap X_{142536}(u) \cap X_{142536}(v)  
\]
 is transverse and consists of $6$ real points.
\end{thm}

Another interesting feature of this problem is that its Galois group~\cite{Ha79}
is not the full symmetric group $S_6$, but rather the symmetric group $S_3$.
This is the smallest known example of an enumerative geometric problem in the Schubert
calculus whose Galois group is not the full symmetric group, and it is strikingly small.

A final phenomena that we discovered is a Schubert problem which is not transverse,
when it involves flags osculating a rational normal curve.
This  may have  negative repercussions for part of
Varchenko's program concerning the Bethe Ansatz and Fuchsian equations~\cite{MV04}. 
It was quite  unexpected, as 
 Eisenbud and Harris showed that on a Grassmannian, any intersection
 \begin{equation}\label{Eq:genInt}
   X_{w_1}(s_1)\cap X_{w_2}(s_2)\cap \dotsb\cap X_{w_m}(s_m)
 \end{equation}
 has the expected dimension
 $\dim(a)-\sum\ell(w_i)$, if the points $s_1,s_2,\dotsc,s_m\in\P^1$ are
 distinct~\cite[Theorem 2.3]{EH83}. 
 Our example shows that the result of Eisenbud and Harris
 cannot be extended to the flag manifold.

 The manifold of flags of type $(1,3)$ in $\C^5$ has dimension 8.
 Since $\ell(32514)=5$ and $\ell(21435)=2$, there are no flags of type $(1,3)$
 satisfying the trio of Schubert conditions $(32514,\,21435,\,21435)$ imposed
 by three general flags.
 If however the flags osculate a rational normal curve $\gamma$, then the
 intersection is nonempty.

\begin{thm}
 $X_{32514}(u)\cap X_{21435}(s)\cap X_{21435}(t)\neq \emptyset$ 
 for all $s,t,u\in\P^1$.
\end{thm}

\begin{proof}
 We may assume that $u=0$, so that flags in $X^\circ_{32514}(u)$ are
 given by matrices in $\calM_{32514}$.
 Consider the $3\times 5$ matrix in  $\calM_{32514}$ with row vectors $v_1$,
 $v_2$, and $v_3$:
 \[
   \left[ \begin{matrix}0&0&1&\frac{3}{2}(s+t)&6st\\
                        0&1&0&-3st&0\\
                        0&0&0&0&1\end{matrix}\right]\ .
 \]
 Let $\Edot\colon E_1\subset E_3$ be the corresponding flag, which lies in 
 $X^\circ_{32514}(\infty)$.
 Define
\[
   \lambda(s)\ :=\ (s^4,\,-4s^3,\,6s^2,\,-4s,\,1)\ \in\ 
    (\C^5)^*\,,
\]
 and note that $\lambda(s)$ annihilates $\gamma(s)$, $\gamma'(s)$,
 $\gamma''(s)$, and $\gamma'''(s)$, where
 $\gamma(s):=(1,s,s^2,s^3,s^4)$ parametrizes the rational normal
 curve in a neighborhood of $\infty$.  Thus $\lambda(s)$ is the linear form
 annihilating the 4-plane $F_4(s)$ osculating the rational normal curve
 $\gamma$ at the point $\gamma(s)$.
 Since $v_1\cdot \lambda(s)^t=0$, we have $E_1\subset F_4(s)$.
 Finally, 
\[
   2sv_1+v_2+(4s^3-6st)v_3\ =\ (0,\,1,\,2s,\,3s^2,\,4s^3)\ =\ \gamma'(s)\,,
\]
 so $E_3\cap F_2(s)\neq 0$ and thus $\Edot\in X_{21435}(s)$.
 Similarly, $\Edot\in X_{21435}(t)$.
\end{proof}

\providecommand{\bysame}{\leavevmode\hbox to3em{\hrulefill}\thinspace}
\providecommand{\MR}{\relax\ifhmode\unskip\space\fi MR }
\providecommand{\MRhref}[2]{%
  \href{http://www.ams.org/mathscinet-getitem?mr=#1}{#2}
}
\providecommand{\href}[2]{#2}

\end{document}